\documentclass{article}
\usepackage{graphicx} % Required for inserting images
\usepackage{amsmath} % cases
\usepackage{amsmath,amsthm}
\newtheorem{theorem}{Theorem}

\newtheorem{clai}{Claim}[theorem]

\newcommand{\cL}{\mathcal{L}}
\newcommand{\cA}{\mathcal{A}}

\title{A characterization of graphs of diameter two with fewer lines than vertices}
\author{Martín Matamala \\
DIM-CMM, CNRS-IRL 2807, Universidad de Chile, Chile\\}
\date{Octubre 2025}

\begin{document}

\maketitle
\begin{abstract}
In 2008 Chen and Chvátal conjectured that any metric space on $n$ points has at least $n$ lines, unless all the points belong to one line. Chvátal proved in 2014 that this is indeed the case for metric spaces with distances $0$, $1$ and $2$. In this work, we prove that there exists a family of ten graphs such that a metric space defined by a graph of diameter two has fewer lines than points if and only if the associated graph belongs to that family. 
\end{abstract}

\section{Introduction}
In a graph $G=(V,E)$, the line defined by $\{x,y\}$, for $x$ and $y$ distinct vertices of $G$, is the set of all vertices $z$ such that there is a shortest path containing $x,y$ and $z$. It is denoted by $\overline{xy}.$
We denote by $\cL(G)$ the set of all lines of $G$.

In this work we characterize graphs $G$ of diameter two with less lines than vertices.
For a graph of diameter two we can determine the distance between two distinct vertices $x$ and $y$ in terms of the set $\{x,y\}$: when it is an edge of the graph then the distance between $x$ and $y$ is 1; otherwise, it is 2.
Therefore, the line defined by two vertices $x$ and $y$ at distance two (equivalently, such that $xy\notin E$) is the set
$$\overline{xy}=\{x,y\}\cup (N(x)\cap N(y)),$$
and when $x$ and $y$ are adjacent, it is given by
$$\overline{xy}=\{x,y\}\cup (N(x)\setminus N(y))\cup (N(y)\setminus N(x)).$$

A result of Chiniforooshan and Chvátal \cite{chichv} showed that 
every graph of diameter two with $n$ vertices has at least $cn^{4/3}$ lines, for some universal constant $c$, and this bound is tight.
On the other hand, by a result of Chvátal \cite{ch2014}, we know that every graph of diameter two either has at least as many lines as vertices or the entire set of vertices is a line. These works arised as partial answers to the following question proposed by Chen and Chvátal \cite{chenchvatal}: {\it 
Is it true that every metric space on $n$ points contains at least $n$ lines or a line containing all the points?} Currently, {the} previous question is known as the Chen-Chvátal conjecture and its validity is still open.

We start by showing ten graphs of diameter two which have less lines than vertices. Our main result is that these are the only graphs of diameter two with this property. 
Although the analysis of these ten graphs is simple, we present it in detail as it is instructive to gain some intuition about the notion of lines. 

Remember that we only consider graphs of diameter two. The smallest such graph is the path of length two, $K_{1,2}$. This graph has only one line which is the entire set of vertices. We shall call such types of lines, \emph{universal lines}.  The same holds true for the cycle of length four, $K_{2,2}$. For the sake of completeness, we mention that it is known that $K_{2,2}$ and paths are the only graphs which define metric spaces having only one line  \cite{shepperd1956transitivities}, which clearly must be a universal line. 

Among the ten graphs, four have five vertices. Let $V=\{x,a,b,c,d\}$ be  the set of vertices of a graph with five vertices. 
We first consider the complete bipartite graph $K_{2,3}$ whose independent sets are $X_1=\{a,c\}$ and $X_2=\{x,b,d\}$. Then, all 2-sets of $V$ define universal lines with the exception of the three subsets of $X_2$. These latter 2-sets define three lines. Therefore, $K_{2,3}$ has only four lines. Here it is also worth to mention that in \cite{maza} it was proved that $K_{2,2}$ and $K_{2,3}$ are the only bridgeless bipartite graphs having less lines than vertices. 

We now consider the complete multipartite graph $K_{1,2,2}$ whose parts are $\{x,b\},\{a,c\}$ and $\{d\}$. Then, the lines $\overline{xa},\overline{xc}, \overline{ab},\overline{cb}$ are the set $\{x,a,b,c\}$. The lines 
$\overline{ad}$ and $\overline{cd}$ are $\{a,c,d\}$, while the lines 
$\overline{xd}$ and $\overline{bd}$ are $\{x,b,d\}$. Finally, the lines $\overline{xb}$ and $\overline{ac}$ are $V$. Therefore, $K_{1,2,2}$ has only four lines. 

Our next graph is the graph 
$K_{1,2,2}'$ which is obtained from $K_{1,2,2}$ by deleting an edge incident to the vertex of degree four, $d$. If the induced cycles of length four of $K_{1,2,2}'$ are $xabc$ and $xadc$, and the cycles of length three are $abd$ and $cbd$, then we have that  $\overline{xa}=\overline{xc}=\overline{ac}=\{x,a,b,c,d\}$, while $\overline{xb}=\overline{ab}=\overline{cb}=\{x,a,b,c\}$ and $\overline{xd}=\overline{cd}=\overline{ad}=\{x,a,c,d\}$. Moreover, we have that 
$\overline{bd}=\{b,d\}$.
As in the previous cases, this graph has only four lines.

The last graph with five vertices we consider is 
$K_{1,2,2}''$, \emph{the house}, obtained from $K_{1,2,2}$ by deleting two edges incident with the vertex of degree four, in such a way that the vertices of degree two do not form an independent set. 
Let us assume that the cycle of length five of $K_{1,2,2}''$ is $xabdc$ and that the chord is $ac$. Then, we have that  $\overline{xa}=\overline{xb}=\{x,a,b\}$, while $\overline{xc}=\overline{xd}=\{x,c,d\}$. Moreover, we have that 
$\overline{ac}=\overline{bd}=\overline{ad}=\overline{bc}=\{a,b,c,d\}$ and 
$\overline{ab}=\overline{cd}=
\{x,a,b,c,d\}=V$.
Again, as in the previous cases, this graph has four lines.

The four remaining graphs are two graphs of six vertices and two of eight vertices. In order to define them, let $K_p$ be a complete graph with $p$ vertices. We can see that $K_p$ has $\binom{p}{2}$ lines: each 2-set defines a different line.  When joining two copies of $K_p$ by a perfect matching $M$ we get a graph $M_{2p}$ where the line generated by two vertices in one copy of $K_p$ is the same as that generated by the corresponding vertices in the other copy. Moreover, the line defined by the ends of an edge on $M$ is the whole set of vertices. The line defined by non-adjacent vertices $a$ and $b'$ in different copies of $K_p$ is the cycle $aa'b'b$, where $aa',bb'\in M$, and thus it is equal to $\overline{ab}=\overline{a'b'}$.
Therefore, the graph $M_{2p}$ has $\binom{p}{2}+1$ lines, when $p\geq 3$. A similar analysis shows that the complete multipartite graph $K_{2p}'$ with all its parts of size two has also $\binom{p}{2}+1$ lines, when $p\geq 3$. Hence, when $p\geq 3$, these graphs have less than $2p$ lines if and only if $p\in\{3,4\}$. In \cite{mapeza} it was proved that the complete multipartite graphs $K_{1,2,2}, K_6'$ and $K_8'$ are the only locally connected graphs which have less lines than vertices. In \cite{abmaroza} the authors proved that 
the graphs $K_{2,2}$, $K_{2,3}$, $K_{1,2,2}$, $K_{1,2,2}'$, $K_6'$ and  $K_8'$ are the only graphs with less lines than vertices in the class $\mathcal{C}$ of all bridgeless graphs $G$ such that every induced subgraph of $G$ is either a chordal graph, has a cut-vertex or a non-trivial module.

Let $\mathcal{F}$ given by 
$$\mathcal{F}=\{K_{1,2},K_{2,2},K_{2,3},K_{1,2,2},K_{1,2,2}',K_{1,2,2}'',M_6,K_6',M_8,K_8'\}
.$$
Therefore, we have proved the backward implication of the following theorem.

\begin{theorem}\label{t:main}
Let $G$ be a graph of diameter two. It has less lines than vertices if and only if it is isomorphic to one of the graphs in $\mathcal{F}$.
\end{theorem}

\begin{proof}
 We now prove the forward implication. Let $G$ be a graph of diameter two with $n$ vertices and less than $n$  lines.

\begin{clai}\label{cl:ind}
Let $S$ be an independent set in $G$ with at least two vertices. Then, for each $x,y\in S$, 
$$\overline{xy}\cap S=\{x,y\}.$$
\end{clai}
\begin{proof}
    Let $z\in S\setminus \{x,y\}$. Then  $xz,yz\notin E$ and thus no shortest path can contain $x,y$ and $z$. Hence $z\notin \overline{xy}$.
\end{proof}
Given a  vertex $x$ in $G$, we denote by $d(x)$ its degree and by $d_2(x)$ the cardinality of the set $N^2(x)$ given by 
$$N^2(x)=\{u\in V\mid ux\notin E\}.$$ With these definitions the set $V$ is the disjoint union of $\{x\}$, $N(x)$ and $N^2(x)$, and thus $$n=|V|=d(x)+d_2(x)+1,$$
for each $x\in V$.

For a given vertex $x$, we consider the following sets of lines
$$\cL^x_1=\{\overline{xa}\mid a\in N(x)\}$$
and 
$$\cL^x_2=\{\overline{xb}\mid b\in N^2(x)\}.$$

Since $G$ has diameter two we have that, for each $b\in N^2(x)$,  $\overline{xb}\cap N^2(x)=\{b\}$. This implies 
\begin{equation*}\label{e:disttwo}
    |\cL^x_2|=|N^2(x)|=d_2(x).
\end{equation*}
Let $x$ be a vertex maximizing $d_2(x)$. 
\begin{clai}
If $d_2(x)=1$, then $G$ is isomorphic to one of the following complete multipartite graphs:  $K_{1,2}$, $K_{2,2}$, $K_{1,2,2}$, $K_6'$  or $K_8'$.
\end{clai}
\begin{proof} 
If $d_2(x)=1$, then each vertex of $G$ has degree at least $n-2$ and $G$ is not a complete graph. Let $u_1,\ldots,u_q$ be the vertices of degree $n-1$ and let $v_1,v_1',\ldots,v_p,v_p'$ be the vertices of degree $n-2$, where $v_i$ and $v_i'$ are non-adjacent for each $i\in\{1,\ldots,p\}$. Then, $p\geq 1$ and $n=2p+q$.

For $i,j\in \{1,\ldots,q\}$, $i\neq j$, the line $\overline{u_iu_j}$ is the 2-set $\{u_i,u_j\}$, as $u_iu_j\in E$ and any vertex not in $\{u_i,u_j\}$ is adjacent to both $u_i$ and $u_j$. Since $p\geq 1$, these types of lines are not universal. For $l,k\in \{1,\ldots,p\}$, $l\neq k$, the line
$\overline{v_lv_k}$ is the set $\{v_l,v'_l,v_k,v_k'\}$, which is also the lines  $\overline{v_lv'_k}$, $\overline{v'_lv_k}$ and $\overline{v'_lv'_k}$. If $(p,q)\neq (2,0)$, these types of lines are not universal. Moreover,
$\overline{v_kv_k'}=V$ and $\overline{u_iv_l}=\overline{u_iv'_l}=\{u_i,v_l,v_l'\}$. If $(p,q)\neq (1,1)$, these last types of lines are not universal.
Therefore, assuming $(p,q)\notin \{(1,1),(2,0)\}$, $G$ has 
$$p(p-1)/2+q(q-1)/2+pq+1$$ lines.
By replacing $q$ for $n-2p$ and doing some arithmetic operations we get that $G$ has 
$$(n-2p)(n-1)/2+p(p-1)/2+1$$ lines. Then, $(n-2p)/2<1$ which implies that $q=n-2p\leq 1$. If $q=1$, then $G$ has 
$(n-1)/2+p(p-1)/2+1$ lines and then $p(p-1)/2+1\leq (n-1)/2=p$ which implies that $p\leq 2$. Since $(p,q)\neq (1,1)$, we get that $p=2$ and thus $G$ is isomorphic to $K_{1,2,2}$.

When $q=0$, then $n=2p$ and $G$ has $p(p-1)/2+1$ lines. Hence, 
$p(p-1)/2+1\leq n-1=2p-1$. From this we get that $p\leq 4$. Since  $(p,q)\neq (2,0)$ and $(p,q)=(1,0)$ is not possible, $p\in \{3,4\}$ and thus $G$ is isomorphic to $K_6'$ or $K_8'$. 
Finally, if $(p,q)=(1,1)$, $G$ is isomorphic to $K_{1,2}$, and if $(p,q)=(2,0)$, $G$ is isomorphic to $K_{2,2}$.
\end{proof}

In the remaining of the proof we consider $d_2(x)\geq 2.$

\subsubsection*{$\cL^x_1\cap \cL^x_2$ is non-empty}
We first analyze the case when $\cL^x_1\cap \cL^x_2$ is non-empty. Let $a\in N(x)$ and $b\in N^2(x)$ such that $\overline{xa}=\overline{xb}$. 

We know that $\overline{xb}\cap N^2(x)=\{b\}$, for each $b\in N^2(x)$. Then, since $N(a)\cap N^2(x)\subseteq \overline{xa}$, we get that $N(a)\cap N^2(x)=\{b\}$. Thus, $d_2(a)\geq d_2(x)-1$ and $d_2(x)\geq 2$ implies that $d(x)\geq 2$.

\begin{clai}Let $a\in N(x)$ and $b\in N^2(x)$ such that $\overline{xa}=\overline{xb}$. Then the line $\overline{ab}$ does not belong to $\cL^x_2$. \end{clai} \begin{proof} In fact, if $\overline{ab}\in \cL^x_2$, then it is $\overline{xb}$. Hence, for each $b'\in N^2(x)$, $b'\neq b$, we have that $b'\notin \overline{xb}$ and thus $b'\notin \overline{ab}$. Since $ab'\notin E$, $b'$ is not a neighbor of $b$. This implies that $N^2(b)$ contains the set $(N^2(x)\setminus \{b\})\cup \{x\})$. Since $G$ has diameter two, for each $b'\in N^2(x)$, $b'\neq b$, there is a vertex $a'$ which is a common neighbor of $b'$ and $a$. As $a$ has no neighbor in $N^2(x)\setminus\{b\}$ we get that $a'\in N(x)\cap N(a)\cap N(b')$. This implies that $a'\notin \overline{xa}$ and thus $a'$ is not adjacent to $b$, since $\overline{xa}=\overline{xb}$. Therefore, we get that $d_2(b)>d_2(x)$ which is a contradiction with the choice of $x$.
\end{proof}
We consider the set of lines $T$ given by 
$$T=\{\overline{ad}\mid d\in N^2(x)\setminus \{b\}\}.$$ 
Then $T$ has $d_2(x)-1$ lines, since $N(a)\cap N^2(x)=\{b\}$, and no line in $T$ contains $x$.

We have that 
$$d_2(a)=d_2(x)-1+d(x)-|N(x)\cap N(a)|-1$$ and that $$d_2(b)=1+|N(x)\cap N(a)|+|N^2(x)\cap N^2(b)|.$$
From this we get that 
$$d_2(a)+d_2(b)=d_2(x)+d(x)+|N^2(x)\cap N^2(b)|-1,$$
which implies that $d_2(x)\geq d(x)+|N^2(x)\cap N^2(b)|-1.$

As the set of lines $$\cL^x_2\cup T\cup \{\overline{ab}\}$$
has $2d_2(x)$ lines and $G$ has at most $d(x)+d_2(x)=n-1$ lines, we get that $d_2(x)\leq d(x)$
and then $|N^2(x)\cap N^2(b)|\leq 1$. 

\begin{clai}
Let $a\in N(x)$ and $b\in N^2(x)$ such that $\overline{xa}=\overline{xb}$.  Then,     $N^2(x)\cap N^2(b)$ is empty.
\end{clai}

\begin{proof} If $N^2(x)\cap N^2(b)$ is not empty, then we get that $d_2(x)=d(x)$. This shows that the set of lines of $G$ is $\cL(G)=\cL^x_2\cup T\cup \{\overline{ab}\},$ which implies that every line in $G$ either contains $x$ or $a$. This is not possible since for each $b'\in N^2(x)\cap N^2(b)$, the line $\overline{bb'}$ neither contains $x$ nor $a$. 
\end{proof} As  $N^2(x)\cap N^2(b)$ is empty we get that $N^2(x)\subseteq N(b)\cup \{b\}$ which implies that the line $\overline{ab}$ is universal, that is, it contains all vertices. It also forces that $$2d_2(x)\geq d_2(a)+d_2(b)=d_2(x)+d(x)-1.$$
Moreover, each line in $T$ contains $b$.

\begin{clai}
    Let $a\in N(x)$ and $b\in N^2(x)$ such that $\overline{xa}=\overline{xb}$ and  $d_2(x)=2$. Then, $G$ is isomorphic to either $K_{1,2,2}''$ or $M_6$.
\end{clai}
\begin{proof}
Let $N^2(x)=\{b,b'\}$.    Since $ab',xb'\notin E$, by the choice of $x$, we get that $b'a'\in E$, for each $a'\in N(x)\setminus \{a\}$. Since $d_2(x)\geq d(x)-1$, the graph $G$ has five or six vertices.   If $d(x)=3$, then $V=\{x,a,a',a'',b,b'\}$ with $N(x)=\{a,a',a''\}$. Again the choice of $x$ implies that we can assume that $ba'\in E$ and since $\overline{xa}=\overline{xb}$, that $aa'\notin E$. From that we get that $aa''\in E$, by the choice of $x$, and then $a''\notin \overline{xa}$ which implies that $a''b\notin E$. 
If $a'a''\notin E$, then $G$ is isomorphic to $M_6$ and we get the conclusion. Otherwise, we prove that $G$ has at least six lines. In fact, when $a'a''\in E$, the lines
$\overline{a'a''}$ and $\overline{aa''}$ do not belong to $$\cL^x_2\cup T\cup \{\overline{ab}\}=\{\overline{xb},\overline{xb'},\overline{ab'},\overline{ab}\}.$$
On the one hand, they do not contain the vertex $x$ which shows that they do not belong to 
$$\{\overline{xb},\overline{xb'},\overline{ab}\}.$$
 On the other hand, $\overline{ab'}=\{a,a'',b,b'\}$, $b'\notin \overline{a'a''}$ and $b',a'\in \overline{aa''}$ which implies that the lines  $\overline{ab'}$, $\overline{a'a''}$ and $\overline{aa''}$ are distinct.

    If $d(x)=2$, we have that $V=\{x,a,a',b,b'\}$ with $N(x)=\{a,a'\}$ and $N(b')=\{b,a'\}$, by the choice of $x$. If $aa'\in E$, then $ba'\notin E$ and $G$ is isomorphic to $K_{1,2,2}''$ with the set $\{a,a',b,b'\}$ inducing a cycle of length four. Otherwise, $aa'\notin E$ which implies that $ba'\in E$ and then $G$ is also isomorphic to the graph $K_{1,2,2}''$ but this time the set $\{x, a,a',b\}$ induces a cycle of length four.    
\end{proof}
When $d_2(x)\geq 3$  we can obtain a complete description of $G$:
\begin{clai}
Let $a\in N(x)$ and $b\in N^2(x)$ such that $\overline{xa}=\overline{xb}$ and $d_2(x)\geq 3$. Then, the graph $G$ is isomorphic to $M_{n}$.
\end{clai}
\begin{proof}
Let $u,v\in N^2(x)$ with $b\notin \{u,v\}$. Then, $a,x\notin \overline{uv}$, as $au,av,xu,xv\notin E$ and thus we have that  $\overline{uv}\notin \cL^x_2\cup T\cup \{ab\}$.   In this case, the set $\cL^x_2\cup T\cup \{\overline{ab},\overline{uv}\}$
contains $2d_2(x)+1$ lines, which implies that it is $\cL(G)$. It also implies that $d_2(x)\leq d(x)-1$ from which we obtain that $d_2(a)=d_2(b)=d_2(x)=d(x)-1$. From this we get that $d(x)=|N(x)\cap N(a)|+2$ which implies that there is $c\in N(x)$ such that 
$$N(x)\setminus N(a)=\{a,c\}=N(b)\cap N(x).$$
Since for each $a'\in N(x)$, $a'\notin  \{a,c\}$, we have that $a,b\notin \overline{xa'}$,
and this latter line belongs to $\cL^x_2\setminus \{\overline{xb}\}$. Hence, there is $b'\in N^2(x)$, $b'\neq b$ such that $\overline{xa'}=\overline{xb'}$. Therefore, by changing the role of $(a,b)$  by $(a',b')$ we conclude that $$N(x)\setminus N(a')=\{a',c\}=N(b')\cap N(x),$$for each $a'\in N(x)$.  This shows that the sets $(N(x)\cup \{x\})\setminus\{c\}$ and $N^2(x)\cup \{c\}$ induce complete graphs, both of the same size and that these complete graphs are joined by the perfect matching $$\{xc,ab\}\cup \{a'b'\mid a'\in N(x)\setminus\{a,c\}, b'\in N^2(x)\setminus\{b\}, \overline{xa'}=\overline{xb'}\}.$$
\end{proof}
{Hence, $G$ is isomorphic to $M_{2p}$ with $p=d_2(x)+1$.}
As we discussed at the {beginning} of this work, for $p\geq 3$, the graph $M_{2p}$ has $1+\binom{p}{2}$ lines and $2p$ vertices. Thus, $1+p(p-1)/2\leq 2p-1$, which implies that $p\leq 4$. {Since $d_2(x)\geq 3$, we get that $p=4$ and thus $G$ is isomorphic to $M_8$ (which has seven lines).}

\subsubsection*{$\cL^x_1\cap \cL^x_2$ is empty}
We now consider the case when $\cL^x_1\cap \cL^x_2$ is empty.

\begin{clai}\label{cl:atwounionathreenonempty} 
If $\cL^x_1\cap \cL^x_2$ is empty, then $|\cL^x_1|<|N(x)|$ .
\end{clai}
\begin{proof}
We first show that $$|\cL^x_1\cup\cL^x_2|<|N(x)|+|N^2(x)|.$$ Otherwise, $|\cL^x_1\cup\cL^x_2|=|N(x)|+|N^2(x)|=n-1$ and we get that $\cL(G)=\cL^x_1\cup \cL^x_2$. In particular, each line of $G$ contains $x$.

Since $d_2(x)\geq 2$, there are distinct $b,b'\in N^2(x)$. Then the line $\overline{bb'}$ does not contain $x$, a contradiction. 

As the graph $G$ has diameter two we know that  $|\cL^x_2|=|N^2(x)|$. 
Therefore, under the assumption that $\cL^x_1\cap \cL^x_2$ is empty we get that 
$|\cL^x_1|<|N(x)|$.
\end{proof}

The previous claim motivates us to define a relation $R$ on $N(x)$ given  by $$uRv\iff \overline{xu}=\overline{xv}.$$ Clearly, $R$ is an equivalence relation. Let us denote by $[u]$ the equivalence class of $u\in N(x)$ with respect to $R$.

\begin{clai}
    For each $u\in N(x)$, $[u]$ is an independent set which is also a module of $G$.
\end{clai}
\begin{proof}
By definition of $R$, we know that for each $v\in [u]$,  $\overline{xu}=\overline{xv}$. Since $u,v\in N(x)$ this implies that $u$ and $v$ can not be adjacent. Therefore, each equivalence class of $R$ is an independent set.

 To see that $[u]$ is a module, we have to prove that for each $z\in V$,  $zu\in E \iff zu'\in E$, for each $u'\in [u]$.  If $z\in N(x)$, then $zu\in E$ if and only if $z\notin \overline{xu}$. This happens if and only if $z\notin \overline{xu'}$, for each $u'\in [u]$. In turns, this latter fact holds if and only if $zu'\in E$, for each $u'\in [u]$. Similarly, if $z\in N^2(x)$, $zu\in E$ if and only if $z\in \overline{xu}$ if and only if $zu'\in E$, for each $u'\in [u]$. 
\end{proof}

\begin{clai}\label{cl:desline}
    Let $\ell=\overline{vw}$ any line in $G$. Then, for each $u\in N(x)\setminus\{v,w\}$, $[u]\subseteq \ell$ or $[u]\cap \ell$ is empty. Moreover, for $v,w\in N(x)$, 
    when $vw\in E$, we have that $[v]\cup [w]\subseteq \ell$ and when $vw\notin E$ we have that  $([v]\cup [w])\cap \ell=\{v,w\}$.
\end{clai}
\begin{proof} 
When $vw\in E$ we know that for each $u\in N(x)\setminus\{v,w\}$,  $u\in \ell$ if and only if $u\in (N(v)\cup N(w))\setminus (N(v)\cap N(w)).$ Since $[u]$ is a module we get that $[u]\subseteq N(v)\setminus N(w)$ or $[u]\subseteq N(w)\setminus N(v)$. Hence, if $u\in \ell$, then $[u]\subseteq \ell$. For the second statement, since $vw\in E$, we have that $[v]\subseteq N(w)\setminus N(v)$  and $[w]\subseteq N(v)\setminus N(w)$,  because $[v]$ and $[w]$ are independent modules of $G$. Hence $[v]\cup [w]\subseteq \ell$.

When $vw\notin E$ we have that for each $u\in N(x)\setminus\{v,w\}$,  $u\in \ell$ implies that $u\in N(v)\cap N(w)$ which in turns implies that $[u]\subseteq N(v)\cap N(w)$ and thus $[u]\subseteq \ell$, since $[u]$ is a module of $G$. Since $vw\notin E$, we have $[v]\cup [w]$ is an independent set which implies that $([v]\cup [w])\cap \ell=\{v,w\}$.
\end{proof}

We partition $N(x)$ in the following three sets $A_i$, $i\in \{1,2,3\}$.
$$A_1=\{u\in N(x)\mid |[u]|=1\},$$
$$A_2=\{u\in N(x)\mid |[u]|=2\}$$
and
$$A_3=\{u\in N(x)\mid |[u]|\geq 3\}.$$
Then, 
$$|N(x)|=|A_1|+|A_2|+|A_3|.$$

From Claim \ref{cl:atwounionathreenonempty}
we know that $|N(x)|>|\cL^x_1|$ which implies that $A_2\cup A_3$ is non-empty. 

Let 
$$A_2'=\{u\in A_2\mid \exists v\in N(x)\setminus [u], uv\notin E\}$$ and 
$A''_2=A_2\setminus A_2'$.

We shall prove that the set $A''_2$ is non-empty. For this purpose, in the following claims, we define two sets of lines, $\cA_3$ and $\cA_2'$, having at least  {$|A_3|$} and {$|A'_2|$} lines, respectively, such that $\cL^x_1\cup \cL^x_2$, $\cA_3$ and $\cA'_2$ are pairwise disjoint.

\begin{clai}\label{cl:athree}
Let $\cA_3$ given by 
$$\cA_3=\bigcup_{u\in A_3}\cA_u,$$
where 
$$\cA_u=\{\overline{u'u''}\mid u',u''\in [u]\},$$
for $u\in A_3$. Then $\cA_3\cap (\cL^x_1\cup \cL^x_2)$ is empty and $|\cA_3|   \geq |A_3|$.

\end{clai}

\begin{proof}
    From Claim \ref{cl:desline} we know that, for each $u\in A_3$,  each line $\overline{u'u''}$ in $\cA_u$ has exactly two vertices in $[u]$: $u',u''$. Thus, $|\cA_u|=\binom{|[u]|}{2}\geq |[u]|.$
From the same claim we have that each line in $\cL^x_1\cup \cL^x_2$ either contains $[u]$ or has no vertex from $[u]$, for each $u\in A_3$. Hence, $\cA_3\cap (\cL^x_1\cup \cL^x_2)$ is empty. 

Let $u,u',v\in A_3$ with $v\notin [u]$ and $u'\in [u]$. Then, $[v]\cap \overline{uu'}\in \{\emptyset, [v]\}$ which implies that $\cA_u\cap \cA_v$ is empty, and then
$$|\cA_3|=|\bigcup_{u\in A_3}\cA_u|=
\sum_{u\in A_3}\frac{|\cA_u|}{|[u]|}\geq \sum_{u\in A_3}1=|A_3|.$$
\end{proof}

\begin{clai}\label{cl:linetwoprime}
Let $\cA'_2$ be given by 
$$\cA'_2=\bigcup_{u\in A'_2}\bigcup_{v\in N(x)\setminus [u], uv\notin E}\cA_{u,v},$$
where
$$\cA_{u,v}=\{\overline{st}\mid s\in [u], t\in [v]\}.$$ 
Then, $\cA'_2\cap (\cL^x_1\cup \cL^x_2\cup \cA_3)$ is empty and $|\cA'_2|\geq |A'_2|.$
\end{clai}
\begin{proof}
From Claim \ref{cl:desline} we know that for each $u\in A_2'$ and $v\in N(x)\setminus [u]$ such that $uv\notin E$,  each line $\ell=\overline{st}\in \cA_{u,v}$ satisfies $\ell\cap([u]\cup [v])=\{s,t\}$
 and $\ell\cap [w]\in \{\emptyset, [w]\}$, for each $w\in N(x)\setminus ([u]\cup [v])$. Then, 
for $u'\in A'_2\setminus ([u]\cup [v]), v'\in N(x)\setminus [u']$ such that $u'v'\notin E$, the sets $\cA_{u,v}$ and $\cA_{u',v'}$ are disjoint.
The same claim shows that 
$\cA_{u,v}\cap (\cL^x_1\cup \cL^x_2\cup \cA_w)$ is empty, for each $w\in A_3$. In fact, a line in 
$\cA_{u,v}$ has exactly one vertex in $[u]$, while a line in $\cL^x_1\cup \cL^x_2\cup \cA_w$ 
either contains $[u]$ or has no vertex in $[u]$. 

We also know that $|\cA_{u,v}|=|[u]||[v]|$ from which we get that $\cA'_2$ 
contains at least $|A'_2|$ lines.
\end{proof}

\begin{clai} If $\cL^x_1\cap \cL^x_2$ is empty, then $N(x)= A''_2$. 
\end{clai}
\begin{proof} We first prove that the set $A''_2$ is non-empty. Let $a_3$ be the number of lines $\overline{xu}$ in $\cL^x_1$ with $u\in A_3$. Then, the set $\cL^x_1$ contains $|A_1|+|A_2|/2+a_3$ lines.
From previous claims we have that the set $\cL^x_1\cup \cL^x_2\cup \cA_3\cup \cA'_2$ contains
at least 
$$|A_1|+|A_2|/2+a_3+d_2(x)+|A_3|+|A'_2|=d(x)+|A_2|/2-|A''_2|+a_3+d_2(x)$$
lines. Since $|\cL(G)|\leq n-1= d(x)+d_2(x)$,  we get that 
$|A_2|/2+a_3\leq |A''_2|$. As $A_2\cup A_3$ is non-empty, we get that $A''_2$ is non-empty.

For the sake of a contradiction, let us assume that the set $N(x)\setminus A''_2$ is non-empty.     
By the definition of $A''_2$ we know that, for each $u\in A''_2$ and $w\in N(x)\setminus [u]$, the edge $uw$ belongs to $E$. Then, for each $v\in A''_2\setminus [u]$ and for each $w\in N(x)\setminus ([u]\cup [v])$, we have that $wu$ and $wv$ belong to $E$. From this and the fact that $[u]$ and $[v]$ are independent sets we get that 
$$\overline{uv}\cap N(x)=[u]\cup [v].$$

Therefore, 
the set of lines 
$$\cA''_2=\{\overline{uv}\mid u,v\in A''_2, u\notin [v]\}$$
has exactly $\binom{|A''_2|/2}{2}$ lines and 
the set of lines 
$$\cA'''_2=\{\overline{uw}\mid w\in N(x)\setminus A''_2, u\in A_2''\}$$
has at least $|A''|/2$ lines, since $N(x)\setminus A''_2$ is non-empty. 
Moreover, $\cA''_2\cap \cA'''_2$ is empty and no line in $\cA''_2\cup \cA'''_2$ contains $x$. 

Therefore the number of lines in the set $$\cL^x_1\cup \cL^x_2\cup \cA_3\cup \cA'_2 \cup \cA''_2\cup \cA'''_2$$
is at least 
\begin{eqnarray*}
    &  &d_2(x)+|A_1|+|A_2|/2+a_3+ |A_3|+|A_2'|+\binom{|A''_2|/2}{2}+|A''_2|/2\\
    &=& d_2(x)+d(x)+|A'_2|/2+a_3+\binom{|A''_2|/2}{2},
\end{eqnarray*}
where the equality comes from the fact that  $$d(x)=|A_1|+|A_2'|+|A_2''|+|A_3| \text{ and } A_2=A_2'\cup A_2''.$$

 Since we have that $d_2(x)+d(x)=n-1\geq |\cL(G)|$ and $A''_2$ is non-empty, we get that $a_3=|A'_2|=0$ and $|A''_2|=2$. Thus  $N(x)=A_1\cup A''_2$, with $A''_2=[u]=[v]=\{u,v\}$ and we also get that 
 $$\cL(G)=\cL^x_1\cup \cL^x_2\cup \cA_3\cup \cA'_2 \cup \cA''_2\cup \cA'''_2=\cL^x_1\cup \cL^x_2\cup \{\overline{uw}\},$$
 with $w\in A_1$.
 
 This shows that for each $b\in N^2(x)$ we have that  $bu,bv\in E$. Otherwise, the line $\overline{bu}$ does not contain $x$ and we have that $\overline{bu}\cap [u]=\{u\} $ which proves that $\overline{bu}$ does not belong to $\cL^x_1\cup \cL^x_2\cup \{\overline{uw}\}$, a contradiction.
 
 From the previous fact, we get that $\overline{uv}=V$. As $\overline{uw}$ does not contain $x$, we get that there is $z\in V\setminus\{x\}$ such that 
 $V=\overline{xz}$. Since
 $u\notin \overline{xw}$ and $w\notin \overline{xu}$, for each $w\in A_1$, we conclude that $z\in N^2(x)$. But this is a contradiction because the fact that $d_2(x)\geq 2$ implies that $V\notin \cL^x_2$.
\end{proof}

\begin{clai}
    If $\cL^x_1\cap \cL^x_2$ is empty and $N(x)=A''_2$, then $A''_2=\{u,v\}$ and $G$ is isomorphic to $K_{2,3}$ or to $K_{1,2,2}'$.
\end{clai}
\begin{proof}

In this case, the set $\cL^x_1\cup \cL^x_2\cup \cA''_2$ has $$|A''_2|/2+d_2(x)+\binom{|A''_2|/2}{2}=
n-1+\binom{|A''_2|/2}{2}-|A''_2|/2
$$ lines, since $d(x)=|A''_2|$. Thus, $\binom{|A''_2|/2}{2}-|A''_2|/2\leq 0$. This  implies that $|A''_2|\in \{2,4,6\}$.

On the one hand, the set of lines $\cL^x_1\cup \cL^x_2\cup \cA''_2$ contains $n-1$ lines, when $|A''_2|=6$, and thus it is $\cL(G)$, and it contains $n-2$ lines, when $|A''_2|=4$. On the other hand, when $|A''_2|=2$ the set $\cA''_2$ is empty and  $\cL^x_1\cup \cL^x_2$ contains $n-2$ lines. Hence, we get that $G$ can have at most one line not in $\cL^x_1\cup \cL^x_2\cup \cA''_2$. This implies that for each $b\in N^2(x)$ and each $u\in A''_2$, the vertices $u$ and $b$ are neighbors. In fact, if $u\notin N(b)$,  then for $v\in [u]$, $v\neq u$, we have that the lines $\overline{ub}$ and $\overline{vb}$  do not contain $x$, from which we get that none of them belongs to $\cL^x_1\cup \cL^x_2$. We also get that $v\notin \overline{ub}$ and $u\notin \overline{vb}$,  which implies that they do not belong to $\cA''_2$ either, since any line in this latter set either contains $[u]$ or does not intersect it. 

Since every vertex in $A''_2$ is adjacent to every vertex in $N^2(x)$, we obtain that each line in $\cL^x_2$ contains $N(x)=A''_2$ and that each line in $\cL^x_1$ contains $N^2(x)$. Moreover, no vertex in $N^2(x)$ belongs to a line in $\cA''_2$ and for $u,v\in A''_2$ with $v\in [u]$ the line $\overline{uv}$ is universal. Since $d_2(x)\geq 2$, no line in $\cL^x_2\cup \cA''_2$ is universal and $\cL^x_1$ contains a universal line only if $|A''_2|=2$. Additionally,  for distinct vertices $b,b'\in N^2(x)$  the line $\overline{bb'}$ does not contain $x$. Therefore, $\overline{uv}\notin \cL^x_2\cup \cA''_2\cup\{\overline{bb'}\}$ and $\overline{bb'}\notin \cL^x_1\cup \cL^x_2\cup \cA''_2$ which shows that  $|A''_2|=2$, since otherwise, the set  $\cL^x_1\cup \cL^x_2\cup \cA''_2$ contains no universal line and then the lines $\overline{uv}$ and $\overline{bb'}$ do not belong to it.  
As $|A''_2|=2$, the set of lines of $G$ is  $\cL^x_2\cup \{\overline{uv},\overline{bb'}\}$. If there are three distinct vertices $b,b'$ and $b''$ in $N^2(x)$, then we must have that $\overline{bb'}=\overline{b'b''}=\overline{b''b}$. Then, we can assume that $b$ and $b'$ are not adjacent and that $b''\in N(b)\cap N(b')$. But then, $u\in \overline{bb'}\setminus \overline{bb''}$. Therefore, $N^2(x)=\{b,b'\}$ and we get that that $G$ is isomorphic to either $K_{2,3}$, when $bb'\notin E$, or to $K_{1,2,2}'$, when $bb'\in E$.

\end{proof}
\end{proof}

\section{Conclusion and future work}
As any path has only one line, there are arbitrarily large graphs with less lines than vertices. In fact, any graph $G$ is contained in a graph with less lines than vertices: it suffices to glue to $G$ a large enough path. However, it is more difficult to find such examples when we forbid the existence of bridges. This leads to the authors in \cite{abmaroza} to conjecture that there are only finitely many bridgeless graphs with less lines than vertices.  

Besides the nine bridgeless graphs presented in this work, we only know three other with less lines than vertices, all of them of diameter three: $M_6'$ and $M_8'$, obtained from $M_6$ and $M_8$ by deleting one edge of the perfect matching between the two complete sets, respectively, and $\hat{M}_8$ obtained from $M_8$ by deleting a perfect matching of one of the complete sets. In a forthcoming work we prove that indeed these are the only \emph{bridgeless} graphs of diameter three with less lines than vertices and we give a complete characterization of the graphs of diameter three with this property; in total they are fourteen and those which are not bridgeless, have vertices of degree one.

\section{Acknowledgements}

This manuscript improves in clarity thanks to a detailed reading of Luciano Villarroel. This work was supported by Basal program FB210005, ANID, Chile.

\bibliographystyle{plain}
\bibliography{bibdiam}

\end{document}